\documentclass[11pt]{article}
\usepackage[a4paper]{anysize}\marginsize{3.5cm}{3.5cm}{1.3cm}{2cm}
\pdfpagewidth=\paperwidth \pdfpageheight=\paperheight
\usepackage{amsfonts,amssymb,amsthm,amsmath,eucal}
\usepackage{pgf}
\usepackage{bbm}
\usepackage{tikz} 
\usepackage{mathrsfs}
\usetikzlibrary{arrows}
\usepackage{color}
\usepackage{float}
\usepackage{subfigure}
\usepackage{caption}

\definecolor{qqqqff}{rgb}{0.,0.,0.}

\theoremstyle{plain}
\newtheorem{thm}{Theorem}[section]
\newtheorem{theorem}[thm]{Theorem}

\newtheorem{lemma}[thm]{Lemma}

\newtheorem{corollary}[thm]{Corollary}
\newtheorem{conjecture}[thm]{Conjecture}

\theoremstyle{definition}
\newtheorem{definition}[thm]{Definition}
\newtheorem{remark}[thm]{Remark}
\newtheorem{example}[thm]{Example}

\newtheorem{thevarthm}[thm]{\varthmname}

\newenvironment{varthm*}[1]{\trivlist\item[]{\bf #1.}\it}{\endtrivlist}

\renewcommand\geq{\geqslant}

\renewcommand\leq{\leqslant}

\newcommand\be{\begin{eqnarray*}}
\newcommand\ee{\end{eqnarray*}}

\newcommand\C{\mathbb C}
\newcommand\Z{\mathbb Z}

\renewcommand\P{\mathbb P}

\newcommand\newop[2]{\def#1{\mathop{\rm #2}\nolimits}}
\newop\log{log}
\newop\ord{ord}
\newop\Gal{Gal}
\newop\SL{SL}
\newop\Bl{Bl}
\newop\mult{mult}
\newop\mass{mass}
\newop\div{div}
\newop\codim{codim}
\newop\sing{sing}
\newop\Zeroes{Zeroes}

\newcommand\alphahat{\widehat{\alpha}}
\newcommand\eps{\varepsilon}

\def\keywordname{{\bfseries Keywords}}%
\def\keywords#1{\par\addvspace\medskipamount{\rightskip=0pt plus1cm
\def\and{\ifhmode\unskip\nobreak\fi\ $\cdot$
}\noindent\keywordname\enspace\ignorespaces#1\par}}
\def\subclassname{{\bfseries Mathematics Subject Classification
(2000)}\enspace}
\def\subclass#1{\par\addvspace\medskipamount{\rightskip=0pt plus1cm
\def\and{\ifhmode\unskip\nobreak\fi\ $\cdot$
}\noindent\subclassname\ignorespaces#1\par}}

\def\endproof{\hspace*{\fill}\endproofsymbol\endtrivlist}

\def\endproofsymbol{\frame{\rule[0pt]{0pt}{6pt}\rule[0pt]{6pt}{0pt}}}

\begin{document}

\author{Grzegorz Malara, Tomasz Szemberg, Justyna Szpond}
\title{On a conjecture of Demailly and new bounds on Waldschmidt constants in $\P^N$}
\date{\today}
\maketitle
\thispagestyle{empty}

\begin{abstract}
   In the present note we prove a conjecture of Demailly for finite sets of sufficiently many
   very general points in projective spaces. This gives a lower bound on Waldschmidt
   constants of such sets. Waldschmidt constants are asymptotic invariants of subschemes
   receiving recently considerable attention \cite{BocFra16}, \cite{BCGHJNSvTT16}, \cite{DHST14}, \cite{FGHLBMS16}, \cite{Har}.
\subclass{MSC 14C20 and MSC 13A15 and MSC 13F20 and 32S25}
\end{abstract}


\section{Introduction}
   In 1980 Jean-Charles Moreau proved the following version of the Schwarz Lemma
   in several complex variables, \cite[Theorem 1.1]{Mor80}.
\begin{theorem}[Moreau]\label{thm:Moreau}
   Let $Z\subset\C^N$ be a finite set of points. For every positive $m\in\Z$,
   there exists a real number $r_m(Z)$ such that for all $R\geq r\geq r_m(Z)$
   and for all holomorphic functions $f$ vanishing to order at least $m$
   at every point of $Z$ there is
   \begin{equation}\label{eq:Moreau}
      |f|_r\leq \left(\frac{2e^{N/2}r}{R}\right)^{\alpha(mZ)}|f|_R,
   \end{equation}
   where $|f|_s=\sup_{|z|\leq s}|f(z)|$ and $\alpha(kW)$ is the least degree
   of a polynomial vanishing at all points of a finite set $W$ to order at least $k$.
\end{theorem}
   The number $\alpha(mZ)$ in the Theorem is optimal, i.e., the statement
   fails with any larger number.
   Several authors, in particular Chudnovsky, were interested in obtaining
   an exponent in \eqref{eq:Moreau} independent of $m$. To this end one
   defines the following quantity \cite{Wal77}.
\begin{definition}[Waldschmidt constant]
   Let $Z\subset\C^N$ be a finite set of points. The \emph{Waldschmidt constant}
   of $Z$ is the real number
   $$\alphahat(Z)=\lim_{m\to\infty}\frac{\alpha(mZ)}{m}.$$
\end{definition}
   The existence of the limit has been showed by Chudnovsky \cite[Lemma 1]{Chu81}.
   It is well known that $\alphahat(Z)=\inf_{m\geq 1}\frac{\alpha(mZ)}{m}.$
   Chudnovsky established also the following fundamental fact, see \cite[Theorem 1]{Chu81}.
\begin{theorem}\label{thm:general bound}
   Let $Z\subset\C^N$ be a finite set of points. Then
   \begin{equation}\label{eq:general bound}
      \alphahat(Z)\geq \frac{\alpha(Z)}{N}.
   \end{equation}
\end{theorem}
   The bound in \eqref{eq:general bound} can now be easily derived from the
   seminal results of Ein, Lazarsfeld and Smith \cite{ELS01}. We discuss it
   briefly below in Section \ref{sec:around}.
   Chudnovsky suspected that the bound in \eqref{eq:general bound} is not
   optimal and raised the following Conjecture, see \cite[Problem 1]{Chu81}.
\begin{conjecture}[Chudnovsky]
   Let $Z\subset\C^N$ be a finite set of points. Then
   \begin{equation}\label{eq:Chudnovsky Conjecture}
      \alphahat(Z)\geq \frac{\alpha(Z)+N-1}{N}.
   \end{equation}
\end{conjecture}
   This has been subsequently generalized by Demailly, see \cite[p. 101]{Dem82}.
\begin{conjecture}[Demailly]\label{conj:Demailly}
   Let $Z\subset\C^N$ be a finite set of points. Then for all $m\geq 1$
   \begin{equation}\label{eq:Demailly Conjecture}
      \alphahat(Z)\geq \frac{\alpha(mZ)+N-1}{m+N-1}.
   \end{equation}
\end{conjecture}
   Of course, for $m=1$ Demailly's Conjecture reduces to that of Chudnovsky.

   There has been recently considerable progress on the Chudnovsky Conjecture
   for general points obtained independently by Dumnicki and Tutaj-Gasi\'nska in \cite{DT16}
   and Fouli, Mantero and Xie in \cite{FMX16}.

   Our main result here is the following.
\begin{varthm*}{Main Theorem}
   The Demailly's Conjecture \eqref{eq:Demailly Conjecture} holds for
   $s\geq (m+1)^N$ very general points in $\P^N$.
\end{varthm*}

\begin{remark}
   For $m=1$ we recover the aforementioned result \cite{DT16} that the Chudnovsky
   Conjecture holds for $s\geq 2^N$ very general points in $\P^N$.
\end{remark}

   Throughout the paper we work over the field $\C$ of complex numbers.

\section{Around the Chudnovsky Conjecture}\label{sec:around}
   Esnault and Viehweg using methods of complex projective geometry have
   proved the following useful result, see \cite[In\'egalit\'e A]{EsnVie83}.
\begin{theorem}[Esnault -- Viehweg]\label{thm:EV}
   Let $I$ be a radical ideal of a finite set of points in $\P^N$
   with $N\geq 2$. Let $m\leq k$ be two integers. Then
   $$\frac{\alpha(I^{(m)})+1}{m+N-1}\leq \frac{\alpha(I^{(k)})}{k},$$
   in particular
   \begin{equation}\label{inq: wc}
   \frac{\alpha(I^{(m)})+1}{m+N-1}\leq \hat\alpha(I).
   \end{equation}
\end{theorem}
   For $N=2$ the inequality in \eqref{inq: wc} establishes Demailly's Conjecture in $\P^2$.
\begin{corollary}\label{cor:Demailly for P2}
   Conjecture \ref{conj:Demailly} holds for arbitrary finite sets of points in $\P^2$.
\end{corollary}
   Around 2000 Ein, Lazarsfeld and Smith established a uniform containment
   result for symbolic and ordinary powers of homogeneous ideals. For the purpose
   of this paper we recall here a somewhat simplified version of their general result.
\begin{definition}[Symbolic power]
   Let $Z=\left\{P_1,\ldots,P_s\right\}$ be a finite set of points in $\P^N$.
   For an algebraic set $W\subset\P^N$, let $I(W)$ be its homogeneous
   defining ideal. Then
   $$I(Z)=I(P_1)\cap\ldots\cap I(P_s)$$
   and for a positive integer $m$
   $$I^{(m)}(Z)=I(P_1)^m\cap\ldots\cap I(P_s)^m$$
   is the $m$th \emph{symbolic power} of $I(Z)$.
\end{definition}
\begin{theorem}[Ein -- Lazarsfeld -- Smith]\label{thm:ELS}
   Let $Z$ be a finite set of points in $\P^N$ and let $I=I(Z)$ be its
   defining ideal. Then the containment
   \begin{equation}\label{eq:ELS}
      I^{(m)}\subset I^r
   \end{equation}
   holds for all $m\geq Nr$.
\end{theorem}
   Theorem \ref{thm:general bound} is an immediate consequence
   of Theorem \ref{thm:ELS}. Indeed, let $Z\subset\C^N\subset\P^N$
   be a finite set of points with the defining ideal $I=I(Z)$. Then
   $$\alpha(I^{(Nr)})\geq \alpha(I^r)=r\alpha(I)$$
   follows from the containment in \eqref{eq:ELS}. Hence
   $$\frac{\alpha(I^{(Nr)})}{Nr}\geq\frac{\alpha(I)}{N}$$
   for all $r\geq 1$. Passing with $r$ to infinity we obtain
   $$\alphahat(I)\geq\frac{\alpha(I)}{N}.$$
\section{A combinatorial inequality}
   In this section we prove the following auxiliary fact.
\begin{lemma}\label{lem:inequality}
   For all $N\geq 3$, $m\geq 1$ and $k\geq m+1$ there is
   $$
      \binom{k(m+N-1)+1}{N}\geq \binom{m+N-1}{N}(k+1)^N.
   $$
\end{lemma}
\proof
   It is convenient to abbreviate $q:=m+N-1$.
   The claim in the Lemma is equivalent to the following inequality
   \begin{equation}\label{eq:product}
      (kq+1)\cdot(kq)\cdot\ldots\cdot(kq+2-N)\geq
      q\cdot(q-1)\cdot\ldots\cdot m\cdot(k+1)^N.
   \end{equation}
   We will group factors in \eqref{eq:product} and show that
   \begin{equation}\label{eq:product i}
   (kq+1-i)(kq+2-N+i)\geq
   (q-i)(m+i)(k+1)^2.
   \end{equation}
   holds for all $i=0,\ldots,\lfloor\frac{N-1}{2}\rfloor$.
   To this end we define
   \begin{equation}\label{eq:function u}
   u(N,m,k,i):=
   (kq+1-i)(kq+2-N+i)-
   (q-i)(m+i)(k+1)^2
   \end{equation}
   and show that this function is non-negative.

   \textbf{Reduction 1.} In the first step, we will show that the difference function
   $$dk(N,m,k,i)=u(N,m,k+1,i)-u(N,m,k,i)$$
   is non-negative.
   Taking this for granted, in order to show that the function in \eqref{eq:function u}
   is non-negative, it suffices to check it for the least allowed value of $k$,
   i.e. for $k=m+1$. In other words the claim in \eqref{eq:function u} reduces
   to the claim that the function
   \begin{equation}\label{eq:function uk}
   uk(N,m,i):=u(N,m,m+1,i)
   \end{equation}
   is non-negative for all $N,m,i$ in the given range.

   Turning to the proof of the Reduction 1 claim, since the difference function
   is linear in $k$, it suffices to show
   \begin{itemize}
   \item[a)] the leading coefficient of $dk(N,m,k,i)$ treated as a polynomial in $k$ is positive and
   \item[b)] the function is non-negative for $k=m+1$.
   \end{itemize}
   The leading coefficient in a) can be written as
   $$(i^2-Ni+i+\frac13N^2)+(\frac23N^2+mN-m-2N+1).$$
   It is elementary to check that the terms in brackets are non-negative.

   Evaluating $dk(N,m,m+1,i)$ we obtain the following expression
   \begin{equation}\label{eq:1}
   2Nm^2-4m^2+2mN^2-4mN-2imN+2i^2m+2im+4m+2N^2-2N+5i-5iN+5i^2
   \end{equation}
   The term $(2N-4)\cdot m^2$ is positive. The remaining summands in \eqref{eq:1}
   can be rearrange in the following way
   \begin{equation}\label{eq:2}
   (2m+2)\cdot N^2-(4m+2+5i+2im)N+(2im+2i^2m+5i+4m+5i^2).
   \end{equation}
   This is a quadratic function in $N$ whose discriminant
   $$\Delta=4-16m^2-12i^2m^2-16m-8im-36i^2m-20i-15i^2$$
   is negative for all $m$ and $i$ in the allowed range. Thus the expression in \eqref{eq:2} is positive.
   This concludes the proof of Reduction 1.

   We study now the function $uk(N,m,i)$ defined in \eqref{eq:function uk}.
   Our approach is similar. We show in

   \textbf{Reduction 2.} that the difference function
   $$dN(N,m,i)=uk(N+1,m,i)-uk(N,m,i)$$
   is non-negative. This follows immediately from the
   following presentation of this function
   \begin{equation}\label{eq:3}
   dN(N,m,i)=m^3+(\frac{N-1}{2}-i)m^2+(2N-4i-2)m+(3\frac{N-1}{2}m^2-3i+1)
   \end{equation}
   Indeed, all terms in brackets in \eqref{eq:3} are non-negative.

   Hence, it is enough to check that the function in \eqref{eq:function u}
   is non-negative for $N=3$ (and $k=m+1$). But this is immediate since
   $$uk(3,m,i)=(m+3)(m+1)(i-1)^2.$$
   This ends the proof of the Lemma.
\endproof

\section{A proof of the Main Theorem}
   In this section we prove the Main Theorem. First we recall from
   \cite[Theorem 3]{DT16} the following crucial observation.
\begin{theorem}[Lower bound on Waldschmidt constants]\label{thm:lower bound}
   Let $Z$ be a set of $s$ very general points in $\P^N$. Then
   $$ \alphahat(Z)\geq\lfloor \sqrt[N]{s}\rfloor.$$
\end{theorem}
   Turning to the proof of the Main Theorem, let $Z$ be a set of $s\geq(m+1)^N$
   very general points in $\P^N$. Since the result holds in $\P^2$
   by Corollary \ref{cor:Demailly for P2}, we may assume here $N\geq 3$.
   There exists a unique integer $k\geq m+1$ such that
   $$k^N\leq s< (k+1)^N.$$
   By Theorem \ref{thm:lower bound} we have $\alphahat(Z)\geq k$.

   We claim that there exists a form of degree $k(m+N-1)-N+1$ vanishing to order at least $m$
   at every point of $Z$. This follows from the dimension count. Indeed, we need to show that
   $$\binom{k(m+N-1)+1}{N}>\binom{N+m-1}{N}\cdot s$$
   holds. Since $s\leq (k+1)^N-1$, it is enough to show that
   $$\binom{k(m+N-1)+1}{N}\geq \binom{N+m-1}{N}\cdot (k+1)^N$$
   holds. This is exactly the statement of Lemma \ref{lem:inequality}.

   It follows that
   $$\alpha(mZ)\leq k(m+N-1)-N+1.$$
   But then
   $$\frac{\alpha(mZ)+N-1}{m+N-1}\leq k\leq \alphahat(Z)$$
   and we are done.
\endproof
   We conclude this note with examples showing that the inequality in Conjecture \ref{conj:Demailly}
   cannot be improved in general. To this end we recall first the notion of star configurations,
   see \cite{GHM13}.
\begin{definition}[Star configuration of points]\label{def:star config}
   We say that $Z\subset\P^N$ is a \emph{star configuration} of degree $d$
   if $Z$ consists of \textbf{all} intersection points
   of $d\geq N$ \textbf{general} hyperplanes in $\P^N$. By intersection points
   we mean the points which belong to exactly $N$ of given $d$ hyperplanes.
\end{definition}
   The assumption \emph{general} in Definition means that
   any $N$ of $d$ given hyperplanes meet in a single point
   and there is no point belonging to $N+1$ or more hyperplanes.
   In particular a star configuration of degree $d$ consists
   of exactly $\binom{d}{N}$ points.
\begin{example}\rm
   Let $Z\subset\P^N$ be a star configuration of degree $d$.
   Then it is easy to check that for any $k\geq 1$
   $$\alpha((1+kN)Z)=(k+1)d-N+1$$
   and hence
   $$\frac{d}{N}=\alphahat(Z)=\frac{\alpha((1+kN)Z)+N-1}{1+kN+N-1},$$
   so that there is equality in \eqref{eq:Demailly Conjecture}
   for infinitely many values of $m=1+kN$.
\end{example}
   The second example is in a sense more exotic.
\begin{example}\rm
   Let $Z$ be the set of points in $\P^2$ defined by the ideal
   $$I=\langle x(y^3-z^3), y(z^3-x^3), z(x^3-y^3)\rangle.$$
   The $Z$ is the union of points
   \begin{equation*}
   \begin{array}{lll}
      P_1=(1:0:0),         & P_2=(0:1:0),         & P_3=(0:0:1),\\
      P_4=(1:1:1),         & P_5=(1:\eps:\eps^2), & P_6=(1:\eps^2:\eps),\\
      P_7=(\eps:1:1),      & P_8=(1:\eps:1),      & P_9=(1:1:\eps),\\
      P_{10}=(\eps^2:1:1), & P_{11}=(1:\eps^2:1), & P_{12}=(1:1:\eps^2).
   \end{array}
   \end{equation*}
   which together with lines
   \begin{equation*}
   \begin{array}{lll}
   L_1:\; x-y=0,        & L_2:\; y-z=0,        & L_3:\; z-x=0,\\
   L_4:\; x-\eps y=0,   & L_5:\; y-\eps z=0,   & L_6:\; z-\eps x=0,   \\
   L_7:\; x-\eps^2 y=0, & L_8:\; y-\eps^2 z=0, & L_9:\; z-\eps^2 x=0.
   \end{array}
   \end{equation*}
   form a $12_3\,9_4$ configuration, see \cite{DST13}.

   The Waldschmidt constant $\alphahat(Z)=3$ has been computed in passing
   in the proof of Theorem 2.1 in \cite{DHNSST15}. In fact the proof shows
   that
   \begin{equation}\label{eq:4}
      \alpha(3k Z)=9k
   \end{equation}
   for all $k\geq 1$.
   We claim now that
   \begin{equation}\label{eq:5}
      \alpha((3k+2)Z)=9k+8
   \end{equation}
   for all $k\geq 0$. For $k=0$ this can be checked computing $I^{(2)}$ explicitly.
   Clearly \eqref{eq:4} implies
   $$\alpha((3k+2)Z)\leq 9k+8.$$
   Indeed, any partial derivative of a polynomial computing $\alpha(3(k+1)Z)$
   has degree $9k+8$ and the right order of vanishing at $Z$.

   Assume that there is a $k\geq 2$ such that
   $$\alpha((3k+2)Z)\leq 9k+7.$$
   Then there is a divisor $D$ of degree $9k+7$ vanishing to order at least $3k+2$
   at every point $P_i$ of $Z$. Intersecting $D$ with any of the lines
   $L_j$ for $j=1,\ldots,9$, we conclude by Bezout Theorem that $L_j$
   is a component of $D$. Hence there exists a divisor $D'=D-\sum_{j=1}^9L_j$
   of degree $9(k-1)+7$ vanishing to order at least $3(k-1)+2$ at every point
   of $Z$. Repeating this argument $k$ times we get a contradiction
   with $\alpha(2Z)=8$.

   Now, for $m=3k+2$ with $k\geq 1$ we obtain the equality in \eqref{eq:Demailly Conjecture}.
\end{example}
\paragraph*{Acknowledgement.}
   These notes originated from discussions during the VIII Workshop on Algebraic Geometry held
   in Lanckorona in October 2016.
   Research of Szemberg and Szpond was partially supported by National Science Centre, Poland, grant
   2014/15/B/ST1/02197.
   Research of Malara was partially supported by National Science Centre, Poland, grant 2016/21/N/ST1/01491.
   We thank Marcin Dumnicki and Micha\l \ Kapustka for helpful conversations.


\bigskip \small

   Grzegorz Malara,
   Department of Mathematics, Pedagogical University of Cracow,
   Podchor\c a\.zych 2,
   PL-30-084 Krak\'ow, Poland.

\nopagebreak
   \textit{E-mail address:} \texttt{grzegorzmalara@gmail.com}

\bigskip
   Tomasz Szemberg,
   Department of Mathematics, Pedagogical University of Cracow,
   Podchor\c a\.zych 2,
   PL-30-084 Krak\'ow, Poland.

\nopagebreak
   \textit{E-mail address:} \texttt{tomasz.szemberg@gmail.com}

\bigskip
   Justyna Szpond,
   Department of Mathematics, Pedagogical University of Cracow,
   Podchor\c a\.zych 2,
   PL-30-084 Krak\'ow, Poland.

\nopagebreak
   \textit{E-mail address:} \texttt{szpond@gmail.com}


\end{document}